\documentclass{article}

\usepackage{amssymb}
\usepackage{amsmath}
\usepackage{color}
\usepackage{enumerate}
\usepackage{enumitem}
\usepackage{titlesec}
\usepackage[unicode,hidelinks]{hyperref}
\usepackage[color]{changebar}

\begin{document}
\newtheorem{Theorem}{Theorem}[section]
\newtheorem{Proposition}[Theorem]{Proposition}
\newtheorem{Lemma}[Theorem]{Lemma}
\newtheorem{Example}[Theorem]{Example}
\newtheorem{Corollary}[Theorem]{Corollary}
\newtheorem{Fact}[Theorem]{Fact}
\newtheorem{Conjecture}[Theorem]{Conjecture}
\newenvironment{Definition} {\refstepcounter{Theorem} \medskip\noindent
 {\mathbf Definition \arabic{section}.\arabic{Theorem}}\ }{\hfill}
\newenvironment{Remarks} {\refstepcounter{Theorem}
\medskip\noindent {\mathbf Remarks
\arabic{section}.\arabic{Theorem}}\ }{\hfill}

\newenvironment{Exercise} {\medskip\refstepcounter{Theorem}
     \noindent {\mathbf Exercise
\arabic{section}.\arabic{Theorem}}\ } {\hfill}

 \newcommand{\qed}{\hfill \ensuremath{\Box}}
\newenvironment{Proof}{{\noindent \bf Proof\ }}{\hfill\qed}

\newenvironment{claim} {{\smallskip\noindent \mathbf Claim\ }}{\hfill}

\def \blue {\color{blue}}
\def\red{\color{red}}
\def \L {{\cal L}}

\def \B {{\cal B}}
\def \mod {{\rm mod \ }}
\def \iso {\cong}
\def \Lor  {\L_{\rm or}}
\def \Lr {\L_{\rm r}}
\def \Lg {\L_{\rm g}}
\def \I {{\cal I}}
\def \M {{\cal M}}\def\N {{\cal N}}
\def \E {{\cal E}}
\def \Proj {{\mathbb P}}
\def \H {{\cal H}}
\def \x {\times}
\def \Stab {{\rm Stab}}
\def \Z {{\mathbb Z}}
\def \V {{\mathbb V}}
\def \C {{\mathbb C}} \def \Cexp {\C_{\rm exp}}
\def \R {{\mathbb R}}
\def \Q {{\mathbb Q}}\def \K {{\mathbb K}}
\def \F {{\cal F}}\def \A {{\mathbb A}}
 \def \X {{\mathbb X}}
\def \G {{\mathbb G}}
\def\HH {{\mathbb H}}
\def \Nn {{\mathbb N}}\def \Nn {{\mathbb N}}
\def\D {{\mathbb D}}
\def \hat {\widehat}
\def \bar{\overline}
\def \Spec {{\rm Spec}}
\def \bul {$\bullet$\ }
\def\proves {\vdash}
\def \Co {{\cal C}}
\def \ACFp {{\rm ACF}_p}
\def \ACF0 {{\rm ACF}_0}
\def \ee {\prec}
\def \Diag {{\rm Diag}}
\def \Diage {{\rm Diag}_{\rm el}}
\def \DLO {{\rm DLO}}
\def \d {{\rm depth}}
 \def \dist {{\rm dist}}
\def \P {{\cal P}}
\def \ds {\displaystyle}
\def \Fp {{\mathbb F}_p}
\def \acl {{\rm acl}}
\def \dcl {{\rm diff}}
\def \dom {{\rm dom}}
\def \tp {{\rm tp}}
\def \stp {{\rm stp}}
\def \Th  {{\rm Th}}
\def\< {\Lngle}
\def \> {\rangle}
\def \n {\noindent}
\def \minusdot{\hbox{\ {$-$} \kern -.86em\raise .2em \hbox{$\cdot \
$}}}
\def\exp {{\rm exp}}\def\ex {{\rm ex}}
\def \td {{\rm td}\ }
\def \ld {{\rm ld}}
\def \span {{\rm span}}
\def \tilde {\widetilde}
\def \d {\partial}
\def \del {\partial}
\def \cl {{\rm cl}}
\def \acl {{\rm acl}}
\def \cN {{\cal N}}
\def \Qalg {{\Q^{\rm alg}}}
\def \th {^{\rm th}}
\def \deg { {\rm deg} }
\def\hat {\widehat}
\def\li {\L_{\infty,\omega}}
\def\lo {\L_{\omega_1,\omega}}
\def\lk {\L_{\kappa,\omega}}
\def \ee {\prec}\def \bSigma {{\mathbf\Sigma}}
\def \mod {\ {\rm mod\ }}
\def \Tor {{\rm Tor}}
\def \td {{\rm td}}
\def \dim {{\rm dim\ }}
 \def \| {\kern -.3em \restriction \kern -.3em}
 \def \lc{\lceil}
 \def \rc{\rceil}
 \def \SR {{\rm SR}}
 \def\L {{\cal L}}
 \def \a{{\mathbf a}}
 \def \b{{\mathbf b}}
 \def\x{{\mathbf x}}
 \def \< { \langle}
 \def \> { \rangle}
 \def \K {{\cal K}}
  \def \MM {{\mathbb M}} 
  \def \vm {{\mathfrak m}}
  
  \def \Der {{\rm Der}}
  \def\A {{\mathbb A}}
  \def\K{{\mathbb K}}
  \def\alg {{\rm alg}}
  \def \ld{{\rm ldim}} \def \dcf{{\rm DCF}} \def\DCF{{\rm DCF}}
  \font \small=cmsy10
\def \dnf{ \hbox{\ \ {\small\char '152}\kern -.65em \lower .4em \hbox{$\smile$}}}
\def \forks{\hbox{\ \ {\small \char '152} \kern -.7em \lower .1em
\hbox{\it /} \kern -1.1em \lower .4em \hbox{$\smile$}}}
\setlist{topsep=-3pt}

\title{Rigid Differentially Closed Fields}
\author{David Marker\footnote{The author was partially supported by the Fields Institute for Research in the Mathematical Sciences}\\Mathematics, Statistics, and Computer Science\\ University of Illinois Chicago} 
\date{}

\maketitle 

\begin{abstract} {\rm Using ideas from geometric stability theory we construct differentially closed fields of characteristic 0 with no non-trivial automorphisms.}
\end{abstract}

\section{Introduction}  

Our goal is to construct countable differentially closed fields of characteristic 0 ($\DCF_0)$ with no nontrivial automorphisms.  We refer to such fields as {\em rigid}.  This answers a question posed by  Russel Miller.  I will say something about Miller's motivation in my closing remarks.  This may at first seem surprising.  One often, naively, thinks of differentially closed fields should behave like algebraically closed fields where there are alway many automorphisms.  Also, differential closures of proper differential subfields always have non-trivial automorphisms.  We sketch the proof of this using ideas from   Shelah's    proof \cite{shelah} of the uniqueness of prime models  for $\omega$-stable theories (see \cite{marker} \S 6.4 or \cite{tz} 9.2). This is a well-known construction.
  
\begin{Proposition} \label{dcl aut} Let $k$ be a differential field with differential closure $K\supset  k$.
Then there are non-trivial automorphisms of $K/k$.
\end{Proposition}
 \begin{Proof}
First note that if  ${\bf d}\in K^n$ and $k\< {\bf d}\> $ is the differential field generated by ${{\bf d}}$ over $k$.  Then $K$
is a differential closure of $k\< {\bf d}\> $. This follows from the fact that in an $\omega$-stable
theory $\M$ is prime over $A\subset\M$ if and only if $\M$ is atomic over $A$ and there are no uncountable sets of indiscenibles (see \cite{tz} 9.2.1).
 
 Let $a\in K\setminus k$.  Since $K$ is the differential closure of  $k$, $\tp(a/k)$ is isolated by some formula $\phi(v)$ with parameters from $k$.  If $a$ is the only element of $K$ satisfying $\phi$, then $a$ is in dcl$(k)=k$, a contradiction. Thus there is $b\in K$ such that $a\ne b$ and $\phi(b)$.

Since $a$ and $b$ realize the same type over $k$ there is $L\models \DCF_0$ with $k\<  b\>  \subseteq L$ and $\sigma:K\rightarrow L$ an  isomorphism such that $\sigma|k$ is the identity and $\sigma(a)=b$.  

 $K$ is a differential closure of both $k\< a\> $ and $k\< b\> $.  Thus $L$ is 
a differential closure of $k\< b\> $ and, by uniqueness of differential closures, there is 
an isomorphism $\tau: L\rightarrow K$ that is the identity on $k\< b\> $.
Then $\tau\circ\sigma$ is an automorphism of $K$ sending $a$ to $b$.
\end{Proof}

\medskip\n {\bf Remarks}\begin{itemize}

\item This argument really shows that if $T$ is an $\omega$-stable theory, 
$A$ is a definably closed substructure of a model of $T$ that is not a model of $T$
and $\M$ is a prime model extension of $A$, then there is a non-trivial automorphism 
of $\M$ fixing $A$ pointwise.

\item While this argument guarantees the existence of a non-trivial automorphism of $K/k$,
it is  possible that it is only one.  If $k$ is a model of Singer's theory of {\em closed ordered differential fields} \cite{singer}, then $k^\dcl=k(i)$ and complex conjugation is the only non-trivial
automorphism of $k^\dcl/k$.  

Omar Le\'on S\'anchez pointed out that the construction of 
a rigid differentially closed fields gives the first known examples of  differentially closed fields $K$ such that $K\ne k(i)$
for any closed ordered differential field $k\subset K$.
 
\item
Proposition \ref{dcl aut} tells us that the rigid differentially closed fields we construct are not the differential closure of any proper differential subfield.

   \end{itemize}

 \medskip
Our construction of rigid differentially closed fields uses ideas from geometric stability theory and work on strongly minimal sets in differentially closed fields of Rosenlicht \cite{rosen}
and Hrushovski and Sokolovi\'c \cite{hs}.  We describe the results we need   in \S 2 and 
construct rigid differentially closed fields in \S3.  We begin \S3 with a warm up constructing   arbitrarily  large rigid models and then give the more subtle construction of rigid countable models.
We refer the reader to \cite{pillay gst} for unexplained model theoretic concepts.  

 I am grateful to Russell Miller for bringing this question to my attention and to Zo\'e Chatzidakis, Jim Freitag, Omar Le\'on S\'anchez and the referee   for remarks on earlier drafts.
 
I am pleased to submit
this paper in honor of Udi Hrushovski's belated 60th birthday.  The main result relies heavily on  his work.

\section{Preliminaries}

We work in $\K\models \dcf$ a monster model of the theory of differentially closed fields of characteristic zero with a single derivation. The constant field $C$ is $\{x\in\K: x^\prime=0\}$.
If $k$ is a differential field and $X\subset \K^n$
is definable over $k$, we let $X(k)$ denote the $k$-points of $X$, i.e., $X(k)=k^n\cap X$.
Of course, by quantifier elimination, $X$ is quantifier free definable over $k$ 

Our main tool will be the strongly minimal sets known as {\em Manin kernels} of elliptic curves.  Manin kernels arose in Manin's proof \cite {manin} of the Mordell Conjecture for function fields in characteristic zero and were central to both Buium's \cite{buium} and Hrushovski's \cite{udi ml} proofs of the Mordell--Lang Conjecture for function fields in characteristic zero.  The model theoretic importance of Manin kernels was developed in the beautiful unpublished preprint of Hrushovski and Sokolovi\'c \cite{hs}.  Proofs of the results from \cite{hs} that we will need all appear in Pillay's survey \cite{pillay dcf} and \cite{marker manin} is another survey on the construction and some of the basic properties of Manin kernels. 

For $a\in K$, let $E_a$ be the elliptic curve $Y^2=X(X-1)(X-a)$.  Let $E_a^\sharp$ be 
the minimal definable differential subgroup of $E$. $E_a^\sharp$ is the closure of 
Tor$(E_a)$ in the Kolchin topology.  

\begin{Theorem}[Hrushovski--Sokolovi\'c] i) If $a^\prime\ne 0$, then $E_a^\sharp$ is a non-trivial locally modular strongly minimal set.

ii) The Manin kernels $E_a^\sharp$ and $E_b^\sharp$ are non-orthogonal if and only if   $E_a$ and $E_b$ are isogenous.  In particular, if $a$ and $b$ are algebraically independent over $\Q$ then 
$E_a^\sharp$ and $E_b^\sharp$ are orthogonal.
\end{Theorem}

In particular, Manin kernels are orthogonal to the field of constants $C=\{x: x^\prime=0\}$.

\medskip

 More generally, if $A$ is a simple abelian variety that is not isomorphic to an abelian variety   defined over the constants we can construct a Manin kernel $A^\sharp$ which is the Kolchin closure of the torsion of $A$ and a minimal infinite definable subgroup of $A$.  $A^\sharp$ is non-trivial locally modular strongly minimal and Hrushovski and Sokolovi\'c also showed that if $X$ is any non-trivial locally modular strongly minimal subset of a differentially closed field, then $X\not\perp A^\sharp$
for some abelian variety $A$.


The other building blocks of our construction are strongly minimal sets introduced by Rosenlicht \cite{rosen} in his proof that the differential closure of a differential field $k$ need not be minimal.

Let $f(X)={X\over 1+X}$.  For $a\ne 0$, let $X_a=\{x: x^\prime=a f(x), x\ne 0\}$.

\begin{Theorem}[Rosenlicht]\label{rosen thm} i) If $a\in k$ and $x\in X_a\setminus k$, then $C(k)=C(k\< x\> )$.

ii)  Suppose $k\subset K$ are differential fields, with $C(K)\subseteq C(k)^\alg$.  Suppose $a,b\in k^\times$, $x\in X_a(K)$, $y\in X_b(K)$ and $x$ and $y$ are algebraically dependent over $k$, then $x,y$ are algebraic over $k$ or $x=y$.
In particular, if $a\ne b$, then $X_a$ and $X_b$ are orthogonal. 
\end{Theorem}

Part i) follows from Proposition 2 of \cite {rosen} while ii) is a slight generalization of Proposition 1 of \cite{rosen} and Gramain \cite{gram}.   These results appear as Theorems 6.12 and 6.2   of \cite{marker dcf}.

\begin{Corollary} Each $X_a$ is a trivial strongly minimal set.
\end{Corollary}\begin{Proof}
By Theorem \ref{rosen thm} i) $X_a$ is orthogonal to the constants.  If $X_a$ were non-trivial, then  $X_a\not\perp A^\sharp$ the Manin kernel of a simple abelian variety.  But if $x\in X_a\setminus k^\alg$, then $k\< x\> =k(x)$ is a transcenence degree 1 extension.
While by results of Buium \cite{buium},  Manin kernels, or anything non-orthogonal to one, give rise
to extensions of transcendence degree at least 2.  Thus $X_a$ is trivial.
\end{Proof}

 \section{Constructing rigid differentially closed fields}
\subsection*{Warm up}

\begin{Proposition}\label{warm up} There are arbitrarily large rigid differentially closed fields.\end{Proposition} 

For this construction
we only need Rosenlicht strongly minimal sets.  
Let $\kappa$ be a cardinal with $\kappa=\aleph_\kappa$.  We will construct a differentially
closed field $K$ of cardinality $\kappa$ such that $|X_a(K)|\ne |X_b(K)|$ for each nonzero $a\ne b$, guaranteeing there is no automorphism sending $a\mapsto b$.

We build a chain of differentially closed fields $K_0\subset K_1\subset\dots\subset K_\alpha\subset\dots$ for $\alpha<\kappa$ 
such that $|K_\alpha|=\aleph_\alpha$.  We simultaneous build $a_0,a_1,\dots,a_\alpha,\dots$
an injective enumeration of $K^\times$ where $K=\bigcup K_\alpha$.  
 
We construct $K$ as follows.

i) $K_0=\Q^\dcl$.  

ii) Given $K_\alpha$ and $a_\alpha\in K_\alpha$. Build $K_{\alpha+1}$ by adding $\aleph_{\alpha+1}$ new independent elements of $X_{a_\alpha}$ and taking the differential closure.

iii) If $\alpha$ is a limit ordinal, let  $K_\alpha=\bigcup_{\beta<\alpha} K_\beta$.\footnote{To build the desired enumeration--let $a_0,a_1,\dots$ be an injective enumeration of $K_0$ and, at stage
$\alpha+1$, let $(a_\gamma: \omega_\alpha\le\gamma<\omega_{\alpha+1})$ be an injective enumeration of $K_{\alpha+1}\setminus K_\alpha$.}

\medskip Since $X_{a_\alpha}\perp X_{a_\beta}$ for $\alpha<\beta$, adding new elements to $X_{a_\beta}$ and taking the differential closure adds no new elements to $X_{a_\alpha}$. Thus $X_{a_\alpha}(K)=X_{a_\alpha}(K_{\alpha+1})$.  In particular $|X_{a_\alpha}(K)|=\aleph_{\alpha+1}$.
Thus there is no automorphism of $K$ with $a_\alpha\mapsto a_\beta$ for $\alpha\ne\beta$.

\medskip One might worry that we have contradicted Proposition \ref{dcl aut}.  Let $B_\alpha$ be
all of the independent realizations  of $X_{a_\alpha}$ that we added at stage $\alpha$.
Then $K$ is the differential closure of $k=\Q\< B_\alpha:\alpha<\kappa\> $.
But, if $b\in X_{a_\alpha}$, then $a_\alpha= {b^\prime (b+1)\over b}\in \Q\< b\> $. Thus $k=K$. 

\subsection*{The countable case}
 
To construct a countable differentially closed field with no automorphisms we will need  
a more subtle mixture of  Rosenlicht extensions with extensions of Manin kernels.
 
Suppose $b\not\in C$.
Let $\dim E_b^\sharp(k)$ be the number of independent realizations in $k$ of the generic type of $E_b^\sharp$ over $\Q\< b\> $.  Manin kernels are useful to us as they can have any countable dimension.  
We will build a countable $K\models \DCF_0$ such that for each $a\ne 0$, there is a natural number 
$$n_a=\max_{b\in X_a(K)} \dim E_b^\sharp(K)$$ and such that $n_a\ne n_b$ for $a\ne b$.
This will guarantee that there is no automorphism with $a\mapsto b$.

Freitag and Scanlon \cite{fs}, and more generally, Casale, Freitag and Nagloo \cite{cfn}, have given constructions of trivial strongly minimal sets which 
can take on any countable dimension.  Presumably these could be used in an alternative construction.

\medskip

We will build $K_0\subset K_1\subset \dots\subset 
K_n\subset\dots$, $a_0,a_1,\dots$ an injective enumeration of   $K^\times=\bigcup K^\times_n$ and $0=n_0<n_1<\dots$ a sequence of natural numbers such
that:
\begin{enumerate}\item\label{0} $C(K_i)=C(K_0)$;
 \item \label{1}$X_{a_i}(K)=X_{a_i}(K_{i+1})$;
\item\label{2} if $b\in X_{a_i}(K)$, then $E_b^\sharp(K)=E_b^\sharp(K_{i+1})$;
\item\label{3} $n_{i+1}= \max_{b\in X_{a_i}(K)} \dim E_b^\sharp$.\footnote{Building the enumeration takes a bit more bookkeeping  in this case.  Let $d_{0,0}, d_{0,1},\dots$ be an injective enumeration of $K_0$ and 
let $d_{i,0}, d_{i,1},\dots$ be an injective enumeration of $K_i\setminus K_{i-1}$.
Start our enumeration of $K$ by letting $a_0=d_{0,0}$.  Suppose we start stage $i$ with 
the partial enumeration $a_0,\dots a_M$. Then for $j=0,\dots i$, let $a_{M+j+1}=d(i, i-j)$.} 
\end{enumerate}
 
 \medskip
If we can do that we will have guaranteed that there are no automorphisms of $K$.

\medskip Let $K_0=\Q^\dcl$.  At stage $s$ we choose a new $a_s\in K_s$.
Let $b_s$ be a element of $X_{a_s}$ generic over $K_s$ and let ${{\bf x}}$ be $n_{s-1}+1$
independent realizations of the generic of $E_{b_s}^\sharp$ over $K_s\< b_s\> = K_s(b_s)$
and let $K_{s+1}=K_s\< b_s,{\bf x}\> ^\dcl$.  

\medskip  By orthogonality considerations, it's clear that conditions 1), 2) and 3) hold, as after 
stage $i+1$ we only add realizations of types orthogonal to $X_{a_i}$ and $E_{b}^\sharp$, for
$b\in X_{a_i}(K)$.  To prove 4) we need to show that there is 
$n_{s}=\max_{d\in X_{a_s}} \dim E_d^\sharp(K_{s+1})$. We have arranged things so that if
there is a bound $n_s$ then $n_s>n_{s-1}$.

\medskip
We need two preliminary lemmas.

\begin{Lemma} \label{dim 0}   If $b^\prime\ne 0$, then  $\dim E_b^\sharp(\Q\< b\> ^\dcl)=0$.   \end{Lemma}
\begin{Proof}
Suppose $x\in E^\sharp_b(\Q\< b\> ^\dcl)$.  All torsion points of $E_b$ are in $\Q(b)^\alg$, so we can suppose  $x$ is a non-torsion point.
But $x$ realizes an isolated type over $\Q \< b\> $.  Let $\psi$ isolate the type
of $x$ over $\Q\< b\> $.  No torsion point can satisfy $\psi$.  Thus by strong minimality
$\psi$ defines a finite set and $x\in \Q\< b\> ^\alg$. \end{Proof}

\medskip Although we will not need it, we can say more in the special case that $\Q\< b\> = \Q(b)$,
for example, if $b\in X_a$ for some $a\in \Q$.  In this case Manin's Theorem of the Kernel \cite{manin} implies that $E_b^\sharp(\Q\< b\> ^\alg)= $ Tor$(E_b)$, see \cite {bp} Corollary K.3.

\begin{Lemma} \label{iso} Suppose $K$ is a differentially closed field, $b,d\in K$ and $E_b$ and $E_d$
are isogenous, then $\dim E_b^\sharp(K)=\dim E_d^\sharp(K)$.
\end{Lemma} 
\begin{Proof} If $E_d$ and $E_b$ are isogenous, then $d$ and $b$ are interalgebraic over $\Q$
and the isogeny $f$ is defined over $\Q(d)^\alg=\Q(b)^\alg$.  Since $f:{\rm Tor}(E_d)\rightarrow
{\rm Tor}(E_b)$ is finite-to-one and the torision is Kolchin dense in a Manin kernel,
$f: E_d^\sharp\rightarrow E_b^\sharp$ is finite-to-one.  It follows that 
$\dim E_d^\sharp(K)=\dim E_b^\sharp(K)$.
\end{Proof}

\medskip
The next lemma shows that we have the necessary bounds.

\begin{Lemma} \label{main}Suppose $K$ is a differentially closed field constructed in a finite iteration
$\Q^\dcl=k_0\subset k_1\subset\dots \subset k_m=K$ where either $k_{i+1}=k_i\< a\> ^\dcl$
where $a$ realizes a trivial type over $k_i$ or $k_{i+1}= k_i\< {\bf x}_i\> ^\dcl$ where ${\bf x}_i$
are $n_i$ independent realizations of the generic type of a Manin kernel $E_{b_i}^\sharp$ 
where $b_i\in k_i$ and 
$E_{b_i}^\sharp \perp E_{b_j}^\sharp$ for $i\ne j$.   If $d\in  K\setminus C$, then $\dim E_d^\sharp(K)=n_i$ for some $i$.
\end{Lemma}
\begin{Proof} We first argue that this is true for each $E_{b_t}^\sharp$.
Define $l_0\subseteq l_1\subseteq\dots \subseteq l_t$ such that
$l_i=k_i\< b_t\> ^\dcl$. Note that $l_t=k_{t}$

By Lemma \ref{dim 0}, $\dim E_{b_t}^\sharp (l_0)=0$.  As we construct $l_1,\dots,l_t$
we are either doing nothing (if $a_i$ or  ${\bf x}_i\in l_{i-1})$ or adding realizations of types
orthogonal to $E_{b_t}^\sharp$.  Thus $\dim E_{b_t}^\sharp(k_t)=0$
and $\dim E_{b_t}^\sharp(k_{t+1})=n_t$.  Since for $i>t$ all $a_i$ and ${\bf x}_i$ realize
types orthogonal to $E_{b_t}^\sharp$, $\dim E_{b_t}^\sharp (K)= n_t$.

Suppose $d\in K\setminus C$.  If $E_d$ is isogenous to some $E_{b_i}$, then, by Lemma \ref{iso}, $\dim E_d^\sharp (K)=\dim E_{b_i}^\sharp(K)=n_i$. So we may assume $E_d^\sharp \perp
E_{b_i}^\sharp$ for all $i$. We claim that in this case $\dim E_d^\sharp(K)=0$.
 For $i\le m$, let $l_i=k_i\< d\> ^\dcl$. By Lemma \ref{dim 0}, $\dim E_d^\sharp (l_0)=0$. 
 As we continue the construction, as above, at each stage we either do nothing or
 realize types that are orthogonal to $E_d^\sharp$. Thus 
 we add no new elements of $E_d^\sharp$ and $\dim E_d^\sharp(K)=0$.\end{Proof}

\bigskip\n We can interweave a many models construction.  In \cite{hs} the authors noted that Manin kernels could be used to show that DCF$_0$ has eni-dop and concluded that there are $2^{\aleph_0}$
non-isomorphic countable differentially closed fields.  An explicit version of this construction coding graphs into models is used in \cite{mm}.  We can fold that coding into our construction of a rigid model.
 
 \begin{Theorem}\label{many} There are $2^{\aleph_0}$ non-isomorphic countable rigid differentially closed fields.  Each of these fields is not the differential closure of a proper differential subfield.
 \end{Theorem}
 
 Consider $X=X_1(\Q^\dcl)$. This is an infinite set of algebraically independent elements.
Let $G=(X,R)$ be a graph with vertex set $X$ and edge relation $R$.  Let $(\{u_i,v_i\}:i=0,1,\dots$ be an enumeration of two element subsets of $X$.
We modify our construction such that at stage $s$ we also add a generic element of 
$E_{u_i+v_i}^\sharp$ if and only if $(u_i,v_i)\in R$.  We can still apply Lemma \ref{main} and 
our construction will produce a rigid differentially closed field $K$. From $K$ we can recover the graph
in an $\L_{\omega_1,\omega}$-definable way.  Thus non-isomorphic graphs will give rise to non-isomorphic rigid differentially closed fields.

\medskip Similarly, we could interweave   graph coding steps in the proof of Proposition \ref{warm up} and build $2^\kappa$ non-isomorphic rigid differentially closed fields of cardinality $\kappa$
when $\kappa=\aleph_\kappa$.

\section{Remarks and Questions}

In \cite {htmmm} and \cite{htmm} the authors introduce the notion of computable and Borel functors between classes of countable structures.  For example, in Theorem \ref{many}, recovering the
graph from the differentially closed field is a Borel functor from differentially closed fields to graphs.
Miller wondered if there could be invertible functors between these classes.  If there is an invertible functor
$F$ from graphs to differentially closed fields, then the authors show that the corresponding automorphism
groups Aut$(G)$ and Aut$(F(G)$ would be isomorphic.  Miller's original idea was that, since there are rigid graphs, one could show there was no such functor by showing that there are no rigid differentially closed fields.  While our construction shows that this idea does not work, never the less, one can show there is no such functor by looking at possible automorphism groups.  It is easy to construct a countable graph with an automorphism of order $n>2$.  But no differentially closed field can have an automorphism of order $n>2$.
Suppose $K$ is differentially closed and $\sigma$ is an automorphism of order $n>2$. Let $F$ be the fixed field of $\sigma$. Then $K/F$ is an algebraic extension of order $n>2$.  By the Artin--Schreier Theorem, this is impossible for $K$ algebraically closed.

\medskip
\n {\bf Question} 1)  Is there a differentially closed field $K$ where $|$Aut$(K)|=2$? 

If so, is the fixed field a model of CODF?  More generally, if $K$ is a real closed differential field
and $K(i)$ is differentially closed, must $K$ be a model of CODF?

\medskip\n {\bf Question} 2) Are there rigid differentially closed fields of cardinality $\aleph_1$?

The construction of such a model would require a new strategy.  Perhaps it would help to assume the set theoretic principle $\diamondsuit$?  Or the methods of \cite{ss msop}.

\end{document}